\documentclass[axioms,article,accept,moreauthors,pdftex,12pt,a4paper]{mdpi}
\usepackage{amsfonts,  amscd}
\usepackage{amssymb}
\lastpage{x}
\doinum{10.3390/------}
\pubvolume{xx}
\pubyear{2015}
\history{Received: xx / Accepted: xx / Published: xx}

 \theoremstyle{mdpi}
 \newcounter{thm}
 \newtheorem{theorem}[thm]{Theorem}
 
 \newtheorem{corollary}[thm]{Corollary}
 \newtheorem{proposition}[thm]{Proposition}
 \theoremstyle{mdpidefinition}
 \newtheorem{definition}[thm]{Definition}
 \newtheorem{definitions}[thm]{Definitions}

\Title{Pro-Lie Groups: A Survey with Open Problems}

\Author{Karl H. Hofmann $^{1,}$* and Sidney A. Morris $^{1,\dagger}$ }

\address{%
$^{1}$ Fachbereich Mathematik, Technische Universit\"at Darmstadt, \hfill\break Schlossgartenstrasse 7, 64289 Darmstadt, Germany, and\hfill\break Department of Mathematics, Tulane University,\hfill\break
New Orleans, LA 70118, USA\\
$^{2}$ Faculty of Science and Technology, \hfill\break
Federation University Australia,
Victoria 3353, Australia, and 
\hfill\break
School of Engineering and Mathematical Sciences,
\hfill\break
La Trobe University, Bundoora, Victoria 3086, Australia}

\contributed{$^\dagger$ These authors contributed equally to this work.}

\corres[2]{\hfill\break hofmann@mathematik.tu-darmstadt.de and morris.sidney@gmail.com.}

\abstract{A topological group is called a pro-Lie group if it is isomorphic to a
closed subgroup of a product of finite-dimensional real Lie groups. This
class of groups is closed under the formation of arbitrary products and
closed subgroups and forms a complete category. It includes each
finite-dimensional Lie group, each locally
compact group which has a compact quotient group modulo its
identity component and thus, in particular, each compact and
each connected locally compact group; it also  includes all
locally compact abelian groups.
This paper provides an overview of
the structure theory and Lie theory of pro-Lie groups including
results more recent than those in the authors' reference book on pro-Lie
groups.  Significantly, it also includes a review of  the recent
insight that  weakly complete unital algebras provide a natural habitat
for both pro-Lie algebras and pro-Lie groups, indeed for the exponential
function which links the two.
(A topological vector space is weakly complete if it is
isomorphic to a power $\R^X$ of an arbitrary
set of copies of $\R$. This class of real vector spaces is at the basis
of the Lie theory of pro-Lie groups.)
The article also lists 12 open questions connected with
pro-Lie groups.
}


\font\euftext =eufm10 scaled \magstep1

\newcommand\RR{\mathbb{R}}
\newcommand\N{\mathbb N}
\newcommand\R{\mathbb R}
\newcommand\T{\mathbb T}
\newcommand\Z{\mathbb Z}
\newcommand\I{\mathbb I}
\newcommand\ssk{\smallskip}
\newcommand\msk{\medskip}
\newcommand\bsk{\bigskip}

\newcommand\defi{\buildrel\rm def 
\over=}

\newcommand\til{\widetilde}        
\DeclareMathOperator{\Hom}{Hom}

 \DeclareMathOperator{\id}{id}

 \newcommand\h{\hbox{\euftext h}}
 \newcommand\g{\hbox{\euftext g}}
 \newcommand\n{\hbox{\euftext n}}
 \renewcommand\r{\hbox{\euftext r}}
 \newcommand\z{\hbox{\euftext z}}
 \renewcommand\L{\hbox{\euftext L}}
 
 \renewcommand\ij{\hbox{\euftext i}}
 \newcommand\sss{\hbox{\euftext s}}
 \newcommand\cc{\hbox{\euftext c}}
 \renewcommand\ll{\hbox{\euftext l}}

\renewcommand\hat{\widehat}         
   \newcommand{\hats}[1]{\ensuremath{\hat{\hat{#1}}}}

 \newcommand\sdir{\rtimes}
\renewcommand\til{\widetilde}
\renewcommand\.{{\cdot}}

\begin{document}


\section{Introduction}

\noindent In 1900  \textsc{David Hilbert} presented a seminal address to the International
Congress of Mathematicians  in Paris. In this address
he initiated a program by
formulating 23 problems  which  influenced a vast amount 
of research of the $20^{\rm th}$ century. The fifth of these problems
asked whether every locally euclidean topological group admits
a Lie group structure. This motivated an enormous volume of 
work on locally compact groups during the first half of the $20^{\rm th}$
century. It culminated in the work of \textsc{Gleason, Iwasawa, 
Montgomery, Yamabe, and Zippin} yielding 
a positve answer to Hilbert's Fifth Problem
and exposing the structure of almost connected 
locally compact groups. (Recall that a topological 
group $G$ is called {\it almost connected} if the quotient group
$G/G_0$, modulo the connected component $G_0$ of the identity, is compact.
The class of  almost connected groups
includes all compact groups and all connected locally compact
groups.) The advances in the second half of the $20^{\rm th}$ century shed much
light on the structure and representation theory
 of locally compact groups.

Notwithstanding this,
it was recognized that the class
of locally compact groups is obviously deficient in that it 
is not even closed under the formation of arbitrary 
products.
 It was unclear which  
class of topological groups would most appropriately extend the class
of finite-dimensional Lie groups and avoid this defect.

At the beginning of the $21^{\rm st}$ century the authors of this survey
introduced the class of {\it pro-Lie groups}. This class 
contains all finite-dimensional Lie groups, all compact groups,
all  locally compact abelian groups and all almost connected 
locally compact  groups. It is characterized as the class 
of all closed subgroups of arbitrary products of finite-dimensional
Lie groups. 
Obviously it is closed under the formation of arbitrary 
products,  and even the passage to closed subgroups.

This notion of pro-Lie group differs from that used 
in the late $20^{\rm th}$ century where a group $G$ is called a 
pro-Lie group if it is {\it locally compact} and contains
arbitrarily small {\it compact} normal subgroups $N$
such that $G/N$ is a finite-dimensional Lie group.

In order to understand the structure of  pro-Lie groups
in our sense 
we developed  a highly  infinite-dimensional Lie theory
of considerable complexity (see \cite{PROBOOK} and subsequent publications
\cite{hm10}, \cite{hm12}, \cite{hm13}). This Lie theory assigns to
 each pro-Lie group $G$  a
pro-Lie algebra $\g$ and
an exponential function $\exp\colon\g\to G$. This approach
was exploited very successfully in \cite{compbook} for compact groups 
and found to be an eminently useful extension of
the classical theory of real Lie groups of finite dimension.

The  theory of an $n$-dimensional real Lie group
is based on the fact that open subsets of $\RR^n$ have a 
rich differentiable structure that is transported to the 
group allowing a differentiable multiplication and inversion
in the group. It has been an ongoing effort to replace
$\RR^n$ by more general, possibly infinite-dimensional, topological
vector spaces supporting differentiable structures.   
The most advanced such theory is the theory of Lie groups on 
differentiable or smooth manifolds modelled on open subsets
of locally convex  vector spaces 
and their real analysis as used by \textsc{H.~Gl\"ockner} 
and  \textsc{ K.-H.~Neeb} (see \cite{gleei}, \cite{gleeii}, \cite{nee}).
One may justifiably ask how the theory of pro-Lie groups 
and the theory of infinite-dimensional
differentiable Lie groups in their spirit are related. 
It was shown  in
\cite{hofnee} that a pro-Lie group is a smooth Lie group in their sense if and
only if it is locally contractible.

The theory of pro-Lie groups has been described in detail in
the 678 page book \cite{PROBOOK} in 2007 and in  later papers.
An endeavor to summarize that work  would therefore be futile. 
Rather, in this survey we highlight central results and explain
some key open problems.

The purely theoretical foundation of, and motivation for, the
theory of pro-Lie groups, however, must be complemented by
an outlook to  areas in which they emerge naturally 
and by necessity. Section 8 therefore deals with 
the appearance of pro-Lie groups in the so-called character
theory of Hopf algebras which received attention in recent
publications \cite{bog}.


\section{The topology of pro-Lie groups}

\noindent It has become a standard result in the literature 
that every connected locally compact
abelian group is homeomorphic to  $\RR^n\times K$ where $K$
is a compact connected abelian group (see e.g.\ \cite{compbook},
Theorem 7.57). The nonabelian version of this 
(see \cite{montzip}, p.~188f.) says that a connected locally compact group
is homeomorphic to $\RR^n\times K$ where $K$ is a compact connected group.
We shall see in this section that a connected pro-Lie group is
homeomorphic to $\RR^{\bf m}\times K$ where ${\bf m}$ is an arbitrary
cardinal and $K$ is a compact connected group. This is convincing evidence
that connected pro-Lie groups represent a natural extension of 
connected locally compact groups. We also see from this observation
 that the pro-Lie group $\RR^I$, for an arbitrary set $I$, 
is a critical example of a pro-Lie
group. We shall return to this theme repeatedly in this text.

\ssk
We shall give now  a complete description of
the topological structure of an almost connected pro-Lie group. We present it
here because it is perhaps easily understood without
too much background information. A complete proof of the result is far
from elementary or short.

A compact connected group $S$ is said to be
{\it semisimple} if its algebraic commutator subgroup is the whole group.
Let us preface the main result by the remark that we know very explicitly
the structure of  compact connected semisimple groups  from
 \cite{compbook}, Theorem 9.19.
It is also a
fact  that in such a group every element is itself a commutator.

Likewise the structure of a  compact connected abelian group $A$ is 
well understood. Indeed a compact abelian group $A$ is connected if and only
if its Pontryagin dual is a torsion 
free abelian group (see \cite{compbook}, Corollary 8.5). This allows the 
determination of details of the structure of such a group as is expounded in
\cite{compbook}, Chapter 8.

We denote by $\Z(n)$ the $n$-element group $\Z/n\Z$.

\ssk
With this notation and information at hand, one can 
appreciate the power of the following result:

\begin{theorem}\label{1.1} {(\rm \cite{hm12}, Corollary 8.9, p.~381)}
 An almost connected pro-Lie group $G$ contains
a compact connected semisimple subgroup $S$ and a compact
connected abelian subgroup 
$A$ such that for suitable sets $I$ and $J$ the  topological group $G$
is homeomorphic to the topological group  \quad
$ \RR^I \times S \times A \times \Delta,$\quad 
where 
$$\Delta=\begin{cases}
\Z(n), &\text{if } G \text{ has finitely many components},\\
                \Z(2)^J, &\text{otherwise.}\\
\end{cases}$$\end{theorem}
\msk

This result allows several immediate corollaries
which are of interest for the topology of pro-Lie groups

\begin{corollary}\label{1.2}
 The space underlying an almost 
connected pro-Lie group is a Baire space. \end{corollary}

This follows from Theorem \ref{1.1} and Oxtoby's results in \cite{ox}.

\begin{corollary}\label{1.3} Every almost connected pro-Lie group
is homotopy equivalent to a compact group.\end{corollary}

Indeed, $\RR^I$ is homotopy equivalent to a singleton. The algebraic
topology of an almost connected pro-Lie group therefore is that of
a compact group. (Cf.\  \cite{compbook}.)

\begin{corollary}
 An almost connected pro-Lie group is
locally compact if and only if in {\rm Theorem \ref{1.1}} the set
$I$ is finite.\end{corollary}

At the root of Theorem \ref{1.1}
is the following main theorem generalizing 
the Theorem  on p.~188ff.\  in \cite{montzip}.

\begin{theorem}\label{1.5} {\rm (\cite{hm12}, Main Theorem 8.1, p.~379)} 
Every almost connected
pro-Lie group $G$ has a maximal compact subgroup $M$, and any other
compact subgroup is conjugate to a subgroup of  $M$. 
Moreover,  
\begin{itemize}
\item[{\rm(1)}] $G=G_0M$,  
\item[{\rm (2)}] $M_0=G_0\cap M$, and
\item[{\rm (3)}] $M_0$ is a maximal compact subgroup of $G_0$. 
\end{itemize}\end{theorem}

We record that the results of Theorem \ref{1.5} enter into a proof
of  Theorem \ref{1.1} in an essential way. 
In the process of proving Theorem \ref{1.1} one also establishes the following
theorem which is more concise than Theorem \ref{1.1} if one assumes 
the structure theory of compact groups as presented in  \cite{compbook}.

\begin{theorem}\label{1.6} {\rm (\cite{hm12}, Theorem 8.4)}   
Let $G$ be an almost connected pro-Lie group,
and $M$ a maximal compact subgroup of $G$. Then $G$ contains a subspace
$E$ homeomorphic to $\RR^I$, for a set $I$, such that
\quad $ (m,e)\mapsto me: M\times E\to G$\quad
is a homeomorphism. Thus $G$ is homeomorphic to $\RR^I\times M$.\end{theorem}

\section{Pro-Lie groups as projective limits of Lie groups}

 \noindent We have  defined pro-Lie groups as closed subgroups of 
products of finite-dimensional real Lie groups.
In fact, they can be equivalently defined as {\it special} 
closed subgroups of
such products, namely, projective limits of finite-dimensional 
Lie groups. 

At first pass this is surprising and its proof requires some effort.

\begin{definitions}\label{2.1}
 A family $\{f_{jk}:G_j\to G_k| j,k \in J\}$
of morphisms of topological groups is called a {\it projective system}
if $J$ is a  directed set satisfying the 
conditions 

(a) for all $j\in J$,
$f_{jj}$ is the identity map of $G_{jj}$, and 

(b) for $i, j, k\in J$ with
$i\le j\le k$ we have $f_{ik}=f_{ij}\circ f_{jk}$.

\ssk
\noindent Given a projective system of morphisms of topological groups, we
define a closed subgroup $G$ of the product 
$\prod_{j\in J}G_j$ to be the set $\{(g_i)_{i\in J}:
g_j=f_{jk}(g_k)$ for all $j\le k$ in $J\}$.
 This group $G$ is called the {\it
projective limit (of the projective system)}, denoted by
 $\lim_{j\in J} G_j$.  \end{definitions}

When we say in the following that a  subgroup $N$ of a
topological group $G$ is {\it a co-Lie subgroup}, we mean
that $N$ is a normal closed subgroup such that the factor group
$G/N$ is a Lie group.

\begin{theorem}\label{2.2}  {\rm(\cite{gloe}, \cite{PROBOOK}, Theorem 3.39 on p.~161, 
\cite{hm10}, \cite{gadfl})} 
For a topological group $G$, the following conditions are equivalent.
\begin{itemize}
\item[{\rm (1)}] $G$ is complete and every identity neighborhood contains
a co-Lie subgroup.
\item[{\rm (2)}] $G$ is a projective limit of Lie groups.
\item[{\rm (3)}] $G$ is a pro-Lie group.
\end{itemize}
\end{theorem}

This theorem, by the way, explains the choice of the name ``pro-Lie group.''
There is a considerable literature on projective limits of finite discrete
groups, called {\it profinite groups} (see \cite{ribza}, \cite{wi}). 
Also amongst the pro-Lie groups 
there are the {\it prodiscrete} groups, namely,
projective limits of discrete groups or, equivalently,
closed subgroups of products of discrete groups. 
There is not much literature on prodiscrete groups.
We formulate the following open question

\ssk
\noindent \textbf{Question 1.} Is there a satisfactory structure theory 
for nondiscrete prodiscrete groups? More particularly, is there a 
satisfactory structure theory even for abelian nondiscrete prodiscrete
groups?

\ssk

We do know that a quotient group of a pro-Lie group 
modulo a closed subgroup
is a pro-Lie group if the quotient group is complete 
(see \cite{PROBOOK}, Theorem 4.1.(i), p.170).
One may reasonably ask whether, for a pro-Lie group $G$ and a closed
normal subgroup $N$, we have sufficient conditions for $G/N$ to be 
complete and therefore a pro-Lie group.

\begin{theorem}\label{3.5} {\rm(\cite{PROBOOK}, Theorem 4.28, p.~202.)} Let $G$ be a pro-Lie group and $N$ a closed normal subgroup of $G$. Each of the 
following condition suffices for $G/N$ to be a pro-Lie group:
\begin{itemize}
\item[{\rm (i)}] $N$ is almost connected and $G/G_0$ is complete.
\item[{\rm (ii)}] $N$ satisfies the First Axiom of Countability.
\item[{\rm (iii)}]$N$ is locally compact.
\end{itemize}
\end{theorem}

The answer to following question is unknown:

\ssk

\noindent\textbf{Question 2.} Let $G$ be a pro-Lie group with identity component $G_0$. 
 Is $G/G_0$ complete (and therefore prodiscrete)?  

\ssk
 
Indeed this is 
the case when $G$ is abelian: See Theorem \ref{3.3}(iii) below.

\section{Weakly complete vector groups}

\noindent In this section we discuss the previously mentioned connected pro-Lie groups
$\RR^I$, for a set $I$. For infinite sets $I$ these
 are the simplest connected pro-Lie groups outside the class of locally
compact groups. However they will appear 
many times in the Lie theory we shall present. 
In particular, they feature in the structure theory of abelian 
pro-Lie groups.

All vector spaces considered here are understood to be vector spaces
over $\RR$. 
For a vector space $E$, let $\Hom(E,\RR)$ denote the set of all linear
functionals on $E$ with the vector space structure and topology it  inherits
from $\RR^E$, the vector space of all functions $E\to\RR$ with the product
topology.

\begin{proposition}\label{3.1}
 For a  topological vector space $V$
the following conditions are equivalent:
\begin{itemize}
\item[{\rm (1)}] There is a vector space $E$ such that the
topological vector spaces  $\Hom(E,\RR)$
and $V$ are isomorphic as topological groups.
\item[{\rm (2)}] There is some set $I$ such that $V$ is isomorphic 
to $\RR^I$ with the product topology.
\end{itemize}
\end{proposition}

\begin{definition}\label{3.2}  A topological vector space $V$
is called {\it weakly complete} if it satisfies the equivalent conditions of
Proposition \ref{3.1}. \end{definition}

Every weakly complete vector space is a nuclear locally convex space (see \cite{schae}, p.~100,
Corollary to Theorem 7.4, p.~103).
The vector space $E$ in Proposition \ref{3.1} is obtained as the vector space
of all {\it continuous} linear maps $V\to \RR$. 
In fact the category of all weakly
complete vector spaces and continuous linear maps between them is dual to 
the category of all  vector spaces and linear maps between 
them in a fashion
analogous to the Pontryagin duality between compact abelian groups 
and discrete abelian
groups. (See \cite{compbook}, Theorem 7.30,\cite{PROBOOK}, Appendix 2: Weakly
Complete Topological Vector Spaces, pp.~629--650.)

\begin{theorem}\label{3.3} {\rm (\cite{PROBOOK}, Theorem 5.20, p.~230)}
  For any abelian pro-Lie group $G$, there is a weakly complete
vector subgroup $V$ and a closed subgroup $H$ such that
(in additive notation)
\begin{itemize}
\item[{\rm (i)}] $(v,h)\mapsto v+h: V\times H\to G$ is an 
isomorphism of topological groups.
\item[{\rm (ii)}] $H_0$ is a compact connected abelian group and every compact
subgroup of $G$ is contained in $H$.
\item[{\rm (iii)}] $H/H_0\cong G/G_0$ and thus $G/G_0$ is prodiscrete.
\item[{\rm (iv)}] If $G_a$ and $H_a$ are the arc components of the identity of
$G$ and $H$, respectively, then $G_a=V\oplus H_a$.
\end{itemize}
\end{theorem}

We note that in the present circumstances, the positive answer to 
Question 2 expressed in conclusion (iii) follows from the compactness of
$H_0$ via Theorem \ref{3.5}.   

\msk

We have seen that products of pro-Lie groups are pro-Lie groups, and
that closed subgroups of pro-Lie groups are pro-Lie groups. 
As a consequence, projective limits of pro-Lie groups are pro-Lie
groups. 

In Section 7 of \cite{bana}, Banaszczyk introduced nuclear abelian groups.
Since all abelian Lie groups are nuclear and the class of nuclear groups
is closed under projective limits, 

{\it all abelian pro-Lie groups are nuclear.}

\ssk

In these circumstances it is somewhat surprising that quotient groups 
of pro-Lie groups
may fail to be pro-Lie groups. Indeed as we shall see in 
the next Proposition \ref{3.4},
there is a quotient group
of a very simple abelian pro-Lie group, namely, of $\RR^{2^{\aleph_0}}$, which
is an abelian topological group that is not complete and therefore is not
a pro-Lie group. However, this situation is not as concerning 
as it might first appear because every quotient group 
of a pro-Lie group has a completion and 
the completion is a pro-Lie group.
 
\msk 

 We consider
the unit interval $\I=[0,1]$ as representative for the sets of
continuum cardinality $2^{\aleph_0}$. Let $\delta_n\in \Z^{(\N)}$ be
defined by 
$$\delta_n=(e_{mn})_{m\in\N},\quad 
e_{mn}= \begin{cases}
             1 & \text{if }m=n,\\ 
              0 &\text{otherwise.}
               \end{cases}$$ 
Then $B=\{\delta_n:n\in\N\}$ generates the free
abelian group $\Z^{(\N)}$ algebraically.

\begin{proposition}\label{3.4} {\rm (\cite{PROBOOK}, Proposition 5.2, p.~214) } 
The free abelian group $\Z^{(\N)}$ of countably
infinite rank has a nondiscrete topology making it a prodiscrete
group $F$,
 so that the following conditions are satisfied:
\begin{itemize}
\item[{\rm (i)}] $F$ is a countable nondiscrete nonmetrizable
pro-Lie group which therefore
is not a Baire space. 
\item[{\rm (ii)}] $F$ can be considered as a closed subgroup of $V=\RR^\I$ 
with dense $\RR$-linear span in $V$,
so that $V/F$ is an  incomplete group whose completion is a compact 
connected and locally connected but not arcwise connected group.
\item[{\rm (iii)}] Every compact subset of $F$ is contained in a
finite rank subgroup. 
\end{itemize}
\end{proposition}

\ssk

The structure theory results we discussed
permit us to derive results on  the duality of
abelian pro-Lie groups. Recall that this class contains
the class of all locally compact abelian groups properly and is
contained in the class of all abelian nuclear groups.

\msk

 For any abelian topological   group $G$ we let  
$\hat G=\Hom(G,\T)$ denote its dual  with the compact open topology, where $\T=\R/\Z$. 
(See e.g. \cite{compbook}, Chapter 7.) There is a natural morphism of abelian
groups $\eta_G\colon G\to\hats G$ given by $\eta_G(g)(\chi)=\chi(g)$
which may or may not be continuous; information regarding this
issue is to be found for instance in \cite{compbook},  notably 
in Theorem 7.7. We shall call an abelian topological  group
{\it semireflexive}
if $\eta_G\colon G\to\hats G$ is bijective and 
{\it reflexive}
if $\eta_G$ is an isomorphism of topological groups; in the latter case
$G$ is also said {\it to have duality}
(see \cite{compbook}, Definition 7.8).

There is  an example of a nondiscrete but 
prodiscrete abelian torsion group due to Banaszczyk (see \cite{bana}, p.~159,
Example 17.11), which is semireflexive but
not reflexive (see also \cite{PROBOOK}, Chapter 14, p.~595, Example 14.15.
Attention: In line 2 of the text of this example, read $N_\alpha=
\Z(2)^{(\{\nu:\ \nu\ge\alpha\})}$, not $\nu<\alpha$).
Therefore we know that the category of abelian pro-Lie groups is not
self-dual under Pontryagin duality.

\begin{theorem}\label{3.6} {\rm  (\cite{PROBOOK} Theorem 5.36, p.~239)}
Every almost connected 
abelian pro-Lie group
is reflexive, and its character group is a direct sum
of the additive topological group of
a real vector space endowed with the finest locally
convex topology and a discrete abelian group.
Pontryagin duality establishes a contravariant
functorial bijection between the categories
of almost connected abelian pro-Lie groups and
the full subcategory of the category of topological
abelian groups containing all direct sums 
of vector groups with the finest locally convex
topology and discrete abelian groups.\end{theorem}

By this theorem, duality
applies to almost connected abelian Lie groups. In particular
we recall that weakly complete topological  vector spaces have 
a good Pontryagin duality.
By Theorem \ref{3.6} above,  the issue of duality of abelian
pro-Lie groups is reduced to groups with compact identity component.
Amongst this class there are all prodiscrete groups. 
In particular, nothing is known about the  duality of prodiscrete
abelian groups. As all abelian pro-Lie groups are nuclear,
whatever is known on the duality of nuclear abelian groups applies
to pro-Lie groups.

\msk
\noindent{\bf Question 3.} Which abelian pro-Lie groups are reflexive? 

\msk
\noindent{\bf Question 4.} Which abelian   pro-Lie groups are 
strongly reflexive in the sense that
all of its subgroups and Hausdorff quotient groups  are reflexive?

\section{The open mapping theorem}

\noindent We have just dealt with the question whether a quotient group
of a pro-Lie group is a pro-Lie groups and we have seen that the
answer is negative sometimes.
We now deal with the question when a surjective morphism of pro-Lie groups
is a quotient  morphism. Specifically, let
 $f\colon G\to H$ be a surjective morphism of pro-Lie groups.
Does this imply that $f$ is an open mapping? This question is answered
negatively in any first course on topological groups by the example
of $G=\RR_d$, the additive group of reals with the discrete topology
and $H=\RR$ with its natural topology. The quest for sufficient
conditions that would secure the openness is answered by any of
the so-called ``Open Mapping Theorems'' in the classical literature 
in Functional Analysis and in the theory of topological groups.
These impose countability conditions on $G$ and a Baire space
hypothethis on $H$. The latter is provided by such properties as
complete metrizability, local compactness, or the pro-Lie group
property. If $\sigma$-compactness is taken as a countability condition
 on $G$  then
the Baire space property of $H$ will force local compactness upon 
$G/\ker f$. The induced bijection $G/\ker f\to H$ is an isomorphism
if and only if  $f$ is open. So settling the issue for bijective $f$
cannot be much of a restriction, assuming that the properties
envisaged  for $G$ are preserved by passing to quotients.

Whichever way the issue is looked at: The proof of an Open Mapping
Theorem for pro-Lie groups is quite different from all other
Open Mapping Theorems we know. 

\ssk

Once more we encounter the 
class of almost connected pro-Lie groups as that class on
which our methods, notably a Lie theory which we have yet to discuss,
yields, an expected result.

\begin{theorem}\label{4.1}  {\rm(\cite{PROBOOK}, 9.60, p.409f.)\cite{hm10}, \cite{gadfl})} 
 Let $G$ and $H$ be pro-Lie groups and  
 let $f\colon G\to H$
be a surjective continuous homomorphism.
Then $f$ is an open mapping if $G$ is almost connected.
In particular, the natural
bijective morphism $G/\ker f\to H$ is an isomorphism of topological
groups.\end{theorem}

The last conclusion yields another instance of a quotient
group of a pro-Lie group which is again a pro-Lie group, giving
us sufficient condition not included in Theorem \ref{3.5} above, namely,
\ssk
{\it $G$ is almost connected and $N$ is the kernel of
a morphism  onto a pro-Lie group.}

\msk
A further corollary of our Open Mapping Theorem is the
following form of the Second Isomorphism Theorem for pro-Lie groups:

\begin{corollary}\label{4.2}  Assume that $N$ and $H$ are two closed 
almost connected subgroups of a topological group with $N$
being normal,  and assume that
$NH$ is a pro-Lie group. Then $H/(N\cap H)\cong NH/N$.
Moreover,  the natural morphism $\mu\colon N\sdir H\to NH$,
$mu(n,h)=nh$ is a quotient morphism (where $H$ is acting as an 
automorphism group on $N$ via inner automorphisms.\end{corollary}

Noting that $\Z^I$ is a pro-Lie group but is not Polish
unless $I$ is countable (see \cite{PROBOOK}, pp.~235, 236,
notably Lemma 5.30), we ask

\ssk
\noindent{\bf Question 5.}\quad  Is a surjective morphism $f\colon\Z^I\to K$ 
onto a compact group open for every set $I$?

\ssk
The results in Theorem \ref{4.1} and in Chapter 5 of \cite{PROBOOK} 
do not answer this
question. If $I$ is countable, then $\Z^I$ is a
Polish group and an Open Mapping Theorem applies
in this case and gives an affirmative answer.
The Open Mapping Theorem for Pro-Lie Groups
does not apply since $\Z^I$ is never
almost connected for $I$ nonempty. 

\ssk
\noindent{\bf Added in Proof}\quad Saak Gabriyelyan and the second author have recently announced a positive answer to Question 5. 

\section{Lie theory}

\noindent We started our discussion with a presentation of some remarkable
structure theorems on almost connected pro-Lie groups. It is not
surprising that the proofs of such results require some powerful tools.
The crucial tool for a structure theory of pro-Lie groups is their
Lie theory. It is a challenge to explain what we mean by ``Lie Theory.''
Minimally one wants to attach to each pro-Lie group $G$ a Lie algebra
$\L(G)$ with characteristics making it suitable for an exploitation of
its topological algebraic structure for the topological group 
structure of the pro-Lie group $G$. 
Guided by our knowledge of classical Lie group theory we
shall link the group $G$ with
its Lie algebra $\L(G)$ by an exponential
function $\exp\colon\L(G)\to G$ which mediates between Lie algebra
theoretical properties of $\L(G)$ and group theoretical properties of $G$.

As a start, for each $X\in\L(G)$ the function 
$t\mapsto \exp(t.X):\RR\to G$ is a morphism of topological groups.
Tradition calls a morphism $f\colon \RR\to G$ of topological groups
a {\it one-parameter subgroup} of $G$ (admittedly, a misnomer!). 
In classical Lie theory, every one-parameter
subgroup is obtained via the exponential function in this fashion.
In other words, however, one defines a Lie algebra and an  exponential
function, it must establish a bijection $\beta$ from the 
elements of the Lie algebra $\L(G)$ to the set $\Hom(\RR,G)$ 
of  one-parameter subgroups of $G$ so that the following
diagram commutes:

$$\begin{CD}
\L(G)    @>\exp>>  G\\
@V\beta VV        @VV\id_GV\\
\Hom(\RR,G)     @>>f\mapsto f(1)>  G
\end{CD}$$

So, if we are given a topological group $G$, as a first step we may
think of $\L(G)$ as $\Hom(\RR,G)$ with a scalar multiplication
such that $r.X(s)=X(sr)$ for $X\in\L(G)$ and $r, s\in\RR$.
If $G$ has additional structure, such as that of a pro-Lie group,
we will obtain additional structure on $\L(G)$.
So if $G$ is a closed subgroup of a product $P=\prod_{i\in I} G_i$ 
of finite-dimensional Lie groups, then an element 
$X\in \L(P)$ may be identified with an
element $(X_i)_{i\in I}$ of $\prod_{i\in I}\L(G_i)$, and an element
$X$ of this kind is in $\L(G)$ simply if $X(\R)\subseteq G$.
As $G_i$ is a finite-dimensional Lie group,
each $\L(G_i)$ has the structure of a Lie algebra, 
$\L(P)$ is both a weakly complete topological
vector space and a topological Lie algebra. Since $G$ is a closed
subgroup of $P$, it is elementary that $\L(G)$ is closed in $\L(P)$
both topologically and in the sense of Lie algebras.

\begin{definition}\label{5.1}  A {\it pro-Lie algebra $\g$} is a topological
real Lie algebra isomorphic to a closed subalgebra of a product of 
finite-dimensional Lie algebras. \end{definition}

Clearly every pro-Lie algebra has a weakly complete topological 
vector space as the underlying topological vector space.
In complete analogy to Theorem \ref{2.2} we have the following characterisation:

\begin{theorem}\label{5.2} {\rm (\cite{PROBOOK}, pp.~138ff.) }
 For a topological Lie algebra $\g$, 
the following conditions are equivalent:
\begin{itemize}
\item[{\rm (1)}] $\g$ is complete and every neighborhood of $0$ contains
a closed ideal $\ij$ such that the Lie algebra $\g/\ij$ is finite-dimensional.
\item[{\rm (2)}]$\g$ is a projective limit of finite-dimensional Lie algebras.
\item[{\rm (3)}] $\g$ is a pro-Lie algebra.
\end{itemize}
\end{theorem}

One notes that our procedure of identifying 
$\L(G)=\Hom(\RR,G)$ for a pro-Lie group
$G$ with the structure of a pro-Lie algebra yields an exponential function
$\exp\colon \L(G)\to G$ by $\exp X=X(1)$ for $X\colon \RR\to G$ in $\L(G)$.
The implementation of this set-up is secured in \cite{PROBOOK}, summarized
in the following:

\begin{theorem}\label{5.3}   To each pro-Lie group, there is 
uniquely and functorially
assigned a pro-Lie algebra $\L(G)$ together with an exponential function
$\exp_G\colon \L(G)\to G$ such that every one-parameter subgroup of $G$
is of the form  $t\mapsto \exp t.X: \L(G)\to G$ with a unique element $X\in \L(G)$ and that the 
subgroup $\langle\exp\L(G)\rangle$ generated by all one parameter subgroups is dense
in the identity component $G_0$ of the identity in $G$.

To each pro-Lie algebra $\g$ there is uniquely and functorially assigned
a connected pro-Lie group $\Gamma(\g)$ with Lie algebra $\g$ and for
each pro-Lie group $G$ with Lie algebra $\L(G)$ permitting an
isomorphism $f\colon \g\to\L(G)$ of pro-Lie algebras there is a unique
isomorphism of pro-Lie groups $f_\Gamma:\Gamma(\g)\to G$
such that the following diagram is commutative with $\exp$ denoting 
the exponential function of $\Gamma(\g)$:

$$\begin{CD}
\g    @>\exp>>  \Gamma(\g)\\
@VfVV        @VVf_\Gamma V\\
\L(G)    @>>\exp_G>  G
\end{CD}$$

A pro-Lie group $G$ is prodiscrete if and only if it is totally disconnected
if and only if  $\L(G)=\{0\}$. Further $\Gamma(\g)$ is always simply connected.
\end{theorem}
\bsk

Sophus Lie's Third Theorem applies  and is perfectly coded into the
existence of the functor $\Gamma$. The exactness properties of the
functor $\L$ are well understood (see \cite{PROBOOK}, Theorem 4.20, p.~188).
Structural results such as we discussed at the beginning of our survey
are all based on a thorough application of the Lie theory of pro-Lie groups.
Since $\L(G)=\L(G_0)$ from the very definition of $\L(G)$, in the strictest
sense it applies to connected pro-Lie groups, but we saw that the essential
facts reach out to almost connected pro-Lie groups.

Of course, since every locally compact group has an open subgroup which
is an almost connected pro-Lie group,
 this Lie theory applies to {\it all} locally compact groups.

\bsk

In classical Lie theory,
Lie algebras are more directly amenable to  structural analysis
than  Lie groups as they are purely algebraic. 
While pro-Lie algebras are both topological and algebraic, they are
neverthless more easily  analyzed than  pro-Lie groups as well.

Let us look as some characteristic features of pro-Lie algebras.

\begin{definitions}\label{5.4} A pro-Lie algebra $\g$ is called
\begin{itemize}
\item[(i)] {\it reductive} if every closed ideal is 
algebraically and topologically a direct summand, 
\item[(ii)]{\it prosolvable} if every finite-dimensional 
quotient algebra is solvable, 
\item[(iii)] {\it pronilpotent} if every finite-dimensional 
quotient algebra is nilpotent.
\end{itemize}

\noindent The center of $\g$ is denoted $\z(\g)$.
A pro-Lie algebra is called {\it semisimple} if it is reductive and
satisfies $\z(\g)=\{0\}$.
\end{definitions}

\begin{theorem}\label{5.5}  {\rm (\cite{PROBOOK} Theorem 7.27, p.~283)}
  A reductive pro-Lie algebra $\g$ is a product of 
finite-dimensional ideals each of which is either simple or else is isomorphic to 
$\R$.

The algebraic commutator algebra $[\g,\g]$ is closed and a product
of simple ideals. Also, 

$$  \g=\z(\g) \oplus [\g,\g],$$
algebraically and topologically.

\end{theorem}
 
\bsk 

If a maximal prosolvable (respectively, pronilpotent) 
ideal of a pro-Lie algebra $\g$ exists, it is called 
the (solvable) {\it radical} $\r(\g)$ (respectively, the
{\it nilradical} $\n(\g)$). If there is a smallest closed
ideal $\ij$ of $\g$ such that $\g/\ij$ is reductive, then
we call it the {\it coreductive radical} $\n_{\rm cored}(g)$
of $\g$.

\begin{theorem}\label{5.6} {\rm (\cite{PROBOOK}, Chapter 7)} Every pro-Lie algebra
has a radical, a nilradical, and a coreductive radical,
and 
$$ \n_{\rm cored}(g) \subseteq \n(\g) \subseteq \r(\g).$$
Moreover,
$$ \n_{\rm cored}(\g) =\overline{[\g,\g]}\cap \r(\g)=
\overline{[\g,\r(\g)]}.$$

There is a closed semisimple subalgebra $\sss$ such that
$\g=\r(\g)\oplus \sss$, where $(r,s)\mapsto r+s:\r(\g)\times\sss\to \g$
is an isomorphism of weakly complete topological vector spaces 
(Levi-Mal'cev decomposition).
Moreover,
$$\overline{[\g,\g]}=\n_{\rm cored}(\g)\oplus\sss.$$
\end{theorem}

\msk
All of this fine structure can be translated to the group level with
due circumspection and ensuing complications. One can get an idea
of this translation from the process of Lie's Third Theorem. 
Among other things, Theorem \ref{5.3} yields
for each of the pro-Lie algebras $\g$ a pro-Lie group 
$\Gamma(\g)$ with an exponential function  
$\exp\colon \g\to \Gamma(\g)$. 
\ssk
If $\g$ is pronilpotent, then $\exp$
is a homeomorphism. In fact, the Baker-Campbell-Hausdorff series
is summable on the weakly complete topological vector space $\g$
yielding a binary operation $\star$, so that for $\Gamma(\g)$
we may take $(\g,\star)$ and for $\exp$ the identity map.
This applies, in particular, to the coreductive radical 
$\n_{\rm cored}(\g)$, for which $\g/\n_{\rm cored}(\g)$ is
reductive.

\ssk

For reductive $\g$, however, the product structure of $\g$
expressed in Theorem \ref{5.3} carries over to a clean product
structure $$ \Gamma(\g)\cong \RR^I\times \prod_{j\in J} S_j$$ for
a family of simply connected simple real Lie groups $S_j$, 
producing in fact a simply connected reductive group $\Gamma(\g)$.

\ssk
These observations show again how far connected pro-Lie groups reach
outside the domain of locally compact connected groups while
their structure remains close to that which is familiar from
finite-dimensional Lie groups due to a fairly lucid 
topological-algebraic structure of pro-Lie algebras.
We note that for every connected pro-Lie group $G$ with
$\L(G)\cong \g$ one has a morphism $f\colon \Gamma(\g)\to G$
with a prodiscrete kernel and a dense image such that
the following diagram is commutative:

$$\begin{CD}
\g    @>\exp_{\Gamma(\g)}>>  \Gamma(\g)\\
@V\cong VV        @VVf V\\
\L(G)     @>>\exp_G>  G
\end{CD}$$

 In \cite{PROBOOK} it is demonstrated that this tool allows a
structural analysis of $G$.

For instance, the existence of the various radicals 
of a pro-Lie algebra has
its correspondence in respective radicals in any connected 
pro-Lie group. For example, every connected pro-Lie group $G$
has a unique largest normal connected solvable subgroup ${\bf R}(G)$,
called its (solvable) {\it radical}. If $G$ is any topological group
whose identity component $G_0$ is a pro-Lie group, then one 
writes ${\bf R}(G)={\bf R}(G_0)$ and also calls ${\bf R}(G)$ the
{\it radical} of $G$.

The volume of additional details of the theory of pro-Lie algebras 
and connected pro-Lie groups presented in \cite{PROBOOK} is immense.
It cannot be expected that a survey such as this  can do complete
justice to it.

\section{ Later developments}

In this section we report on some developments in the theory of 
pro-Lie groups since the appearence of \cite{PROBOOK}. Of course we have already
included some important material which appeared subsequent to \cite{PROBOOK},
namely, \cite{hm12}  and \cite{hofnee}.

\msk

The article \cite{hm12} contributed the insight that some essential
structure theorems on connected pro-Lie groups could be formulated
so as to include  almost connected pro-Lie groups. This
provided a common generalization of the structure theories of
connected pro-Lie groups and compact groups. This generalization
is both significant and satisfying.
In this survey this was illustrated by  
Theorem \ref{1.1} and its corollaries and Theorems \ref{1.5} and \ref{1.6}. 

\msk
 
\cite{hm12} contains another interesting result which we think has yet
to be exploited in the literature.

\begin{theorem}\label{6.1} Let $G$ be an arbitrary topological group whose 
identity component $G_0$ is a pro-Lie group. Then there is a closed
subgroup $G_1$ whose identity component is the radical ${\bf R}(G)$ 
such that the following conditions hold:\begin{itemize} 
\item[{\rm (i)}] $G=G_0G_1$ and $G_0\cap G_1={\bf R}(G)$
\item[{\rm (ii)}] The factor group $G/{\bf R}(G)$ is the semidirect
product of the connected normal subgroup $G_0/{\bf R}(G)$ and the 
totally disconnected closed subgroup $G_1/{\bf R}(G)$.
\item[{\rm (iii)}] In particular, 
$$\frac G {G_0}\cong\frac{G/{\bf R}(G)}{G_0/{\bf R}(G)}
                  =\frac{G/{\bf R}(G_0)}{G_0/{\bf R}(G_0)}
\cong\frac{G_1}{(G_1)_0}$$
with a prosolvable pro-Lie group $(G_1)_0$. 
\end{itemize}
\end{theorem}

The significance of Theorem \ref{6.1} emerges even when it is specialized
to the case that $G$ is a pro-Lie group. As was emphasized by formulating
Question 2, we do not know whether the component factor group $G/G_0$
is complete and therefore is a prodiscrete group. Theorem \ref{6.1} reduces
the problem to the case that the identity component of $G$ is prosolvable.
For instance, we obtain a positive answer to Question 2 if we know that
the radical {\bf R}$(G)$ is locally compact or first countable (see 
Theorem \ref{3.5} above). 

In the process of extending the structure theory of pro-Lie groups
from connected ones  to almost connected ones, G. Michael, A. A.
has proved the following structure theorem guaranteeing a local
splitting---provided the nilradical is not too big on the Lie algebra
side.

\begin{theorem}\label{6.2} {\rm (\cite{hm13}, \cite{gadfll})} Assume that $G$ is an 
almost connected pro-Lie group $G$ whose Lie algebra $\g$ has a finite-dimensional coreductive radical $\n_{\rm cored}(\g)$. Then there are
arbitrarily small closed normal subgroups  $N$ such that there exists a 
simply connected Lie group $L_N$ and a morphism
$\alpha\colon L_N\to G$ such that the morphism 
$$(n,x)\mapsto n\alpha(x): N \times L_N\to G$$
is open and has a discrete kernel. In particular, $G$ and
$N\times L_N$ are locally isomorphic. 
\end{theorem}

Let us recall Iwasawa's Local Splitting Theorem for locally compact
groups as it is reported in \cite{compbook}, Theorem 10.89.

\begin{theorem}\label{6.3} Let $G$ be a locally compact group. Then there 
exists an open almost connected subgroup $A$ such that for each 
identity neighborhood $U$ there is 
\begin{itemize}
\item[---]a compact normal subgroup $N$ of $A$ contained in $U$,
\item[---] a simply connected Lie group $L$, and
\item[---] an open and continuous surjective morphism 
          $\phi\colon N\times L\to A$ with discrete kernel
          such that $\phi(n,1)=n$ for all $n\in N$.
\end{itemize}
\end{theorem}

The way to compare the two preceding theorems is to look at
the Lie algebra $\g$ of $G$ in Theorem \ref{6.3}.
We notice that $\g=\L(N)\oplus \ll$, $\ll=\L(L)$. 
As the Lie algebra of a compact
group $N$ of the first direct summand is of the form $\L(N)=\cc\oplus \sss_1$ with
a central ideal $\cc$ and a compact semisimple ideal $\sss_1$. 
The finite-dimensional Lie algebra $\ll$ has a Levi-Mal'cev decomposition
$\r\oplus \sss_2$ with its radical $\r$ and a finite-dimensional
semisimple subalgebra $\sss_2$, so that we have 
$$\g=(\cc\oplus\sss_1)\oplus(\r+_{\rm sdir}\sss_2)
    =(\r\oplus \cc) +_{\rm sdir} (\sss_1+\sss_2).$$
We observe that the radical $\r(\g)$ is $\r\oplus \cc$ and that
the coreductive radical $\n_{\rm cored}(\g)$ of $\g$ is contained in
the finite-dimensional subalgebra $\r$ and is, therefore, finite-dimensional. The hypothesis that was imposed in Theorem \ref{6.2},
namely, that the Lie algebra dimension of 
the coreductive radical is finite,  thus
emerges as a {\it necessary} condition in the more classical Iwasawa
Local Decomposition Theorem of Locally Compact Groups.
Additional comments on the local decomposition of pro-Lie groups
may be found in\cite{hm9}.

\section{A natural occurence of pro-Lie groups and pro-Lie algebras}

\noindent
We emphasized that weakly complete topological vector spaces play an
essential role in the theory of pro-Lie groups and pro-Lie algebras.
Now we record that they are crucial in describing a  mathematical
environment, 
where they occur naturally; this was  pointed out recently in \cite{bog}. 
Each of the categories of real vector spaces and of its dual
category of weakly complete topological vector spaces (see Section 3,
Definition 3.2ff.) in fact has a tensor product (see \cite{bog}, Appendix C, 
notably C4, and \cite{dah}) making each of them 
into a {\it commutative monoidal category} (see e.g. \cite{compbook},
Appendix 3, 
Definition A3.62 ff.). Let us denote by $\mathcal{ V}$ 
the category of (real) vector
spaces $E$, $F\dots$, and of linear maps, equipped
 with the usual tensor product $E\otimes F$, and let us
call $\mathcal{ W}$
  the category of weakly complete vector spaces 
$V$, $W,\dots$ etc. with continuous linear maps,
equipped with  the (completed) tensor product 
$V\til\otimes W$. Then a {\it monoid} in $\mathcal{W}$
(see  \cite{compbook} Appendix 3, Discussion preceding Definition A3.64.)
is a weakly complete topological (associative) algebra with identity,
specifically, a morphism $m\colon A\til\otimes A\to A$ plus a morphism
$\RR\to A$ representing the identity (see also \cite{dah}). 
Its dual $E\defi A'=\Hom(A,\RR)$ then
is an (associative unital)
 {\it coalgebra $m'\colon E\to E\otimes E$ with a coidentity
(augmentation) $u\colon E\to \R$}. The theory of coalgebras (see  \cite{wmich})
culminates in  one theorem holding without any further hypotheses:

\begin{theorem}\label{7.1} {\rm (The Fundamental Theorem of Coalgebras,\cite{wmich} 4.12,
p.~742)} Each associative unitary coalgebra is the directed union
of its finite-dimensional unitary subcoalgebras. \end{theorem}

By duality this  implies at once  \textsc{Rafael Dahmen}'s Fundamental Theorem
of Weakly Complete Algebras.

\begin{theorem}\label{7.2} {\rm (Fundamental Theorem of Weakly Complete Algebras,
\cite{bog}) }For every weakly complete associative unital  algebra $A$
there is a projective system of surjective linear
morphisms $f_{jk}\colon A_k\to A_j$, $j\le k$, $j,k\in J$ of finite-dimensional associative unital algebras and a natural isomorphism
 $\phi_A\colon A\to\lim_{j\in J}A_j$ onto  the projective limit 
of this system. \end{theorem}

Conversely, by definition, every projective limit
   of finite-dimensional real unital associative algebras is
   a weakly complete  associative unital algebra.

\ssk
Every associative
algebra $A$ becomes a Lie algebra $A_{\rm Lie}$ when it is equipped with
the commutator bracket $(x,y)\mapsto [x,y]\defi xy-yx$. 
Each of the Lie algebras $(A_j)_{\rm Lie}$ with $j\in J$ is finite-dimensional. 
From 
Theorem \ref{7.2} and the Definition \ref{5.1} plus Theorem \ref{5.2} we thus have

\begin{corollary}\label{7.3} For every weakly complete associative unital
algebra $A$, the Lie algebra $A_{\rm Lie}$ is a pro-Lie algebra and 
$\phi_A\colon A_{\rm Lie}\to\lim_{j\in J}(A_j)_{\rm Lie}$
is an isomorphism of pro-Lie algebras. \end{corollary}

Each of the morphisms $f_{jk}$ maps the group $A_k^{-1}$ of units,
that is, invertible elements, into the  group of units $A_j^{-1}$,
and so we obtain a natural isomorphism 
$$\phi_A|A^{-1}:A^{-1}\to\lim_{j\in J} A_j^{-1}.$$ 
Now every $A_j^{-1}$ is a (linear)  Lie group, (see\cite{PROBOOK}, Chapter 5, 
5.1--5.32) with exponential function 
$\exp_{A_j}:(A_j)_{\rm Lie}\to A_j^{-1}$. 
By Theorem \ref{2.2}, as a consequence we have \textsc{Dahmen}'s corollary

\begin{corollary}\label{7.4} {\rm (\cite{bog}, Proposition 5.4)}
  The group $A^{-1}$ of units of every 
weakly complete associative unital algebra $A$ is 
a pro-Lie group with Lie algebra $A_{\rm Lie}$ and
exponential function
$$\exp_A:A_{\rm Lie}\to A^{-1},\quad \exp(x)=1+x+\frac 1 {2!}x^2
    +\frac 1 {3!}x^3+\cdots.$$
\end{corollary}

What we have exposed here is the basis of a theory that applies 
to {\it group objects} in the commutative monoidal category $\mathcal{W}$
as defined in \cite{compbook}  Definition A3.64(ii). These objects are
commonly called {\it Hopf algebras} and so we shall fix the following
definition.

\begin{definition}\label{7.5} A {\it weakly complete Hopf algebra} is
a group object in $\mathcal{W}$ according to Definition A3.64(ii) 
of \cite{compbook}.\end{definition}

In particular, the multiplicative structure of a weakly complete
Hopf algebra $A$ is a weakly complete associative unital algebra.
It also has a comultiplication $c\colon A\to A\til\otimes A$ 
linked with the multiplication $m\colon A\til\otimes A\to A$ 
through several commutative diagrams for which we 
refer to \cite{compbook}, A3.63(ii) and which express the fact
that $c$ is indeed a morphism of algebras.
 
\begin{definition}\label{7.6} Let $A$ be a weakly complete Hopf algebra
with a comultiplication $c$. An element $x\in A$ is called
{\it grouplike} if 
it satisfies $c(x)=x \til\otimes x$ and $u(x)=1$,
and it is called {\it primitive} if it satisfies 
$c(x)=x\til\otimes 1+ 1\til\otimes x$. The set of
grouplike (respectively, primitive)
elements will be denoted by $G(A)$ (respectively, $P(A)$).\end{definition}

\msk

One shows the following fact:

\begin{theorem}\label{7.7} $G(A)$ is a closed subgroup of $A^{-1}$, the
group of units of the underlying algebra, and $P(A)$ is a
closed Lie subalgebra of $A_{\rm Lie}$. \end{theorem}

The link to the previous remarks is provided by
the following theorem:

\begin{theorem}\label{7.8} Let $A$ be a weakly complete Hopf algebra.
Then the set of grouplike elements $G(A)$ is a pro-Lie group
and the set of primitive elements  $P(A)$ is a pro-Lie algebra
and is the Lie algebra of $G(A)$ with the exponential function
$$ \exp_A|P(A): P(A) \to G(A)$$
with the restriction of the exponential function of 
{\rm Corollary \ref{7.3}}.\end{theorem}

A proof is based on the fact that for the algebra morphism $c$ we have
$\exp_A \circ c =c \circ (\exp_A\til \otimes \exp_A)$.
If $x$ is primitive, then $c(x)=(x\til\otimes1)+(1\til\otimes x)$
and thus $c(\exp x)=\exp(c(x))=\exp((x\til\otimes1)+(1\til\otimes x))
=(\exp x\til\otimes1)(1\til\otimes\exp x)=\exp x\til\otimes \exp x$. 
Therefore $\exp$ maps $P(A)$ into $G(A)$. To see the converse,
let $t\mapsto \exp t\.x\colon \R\to G(A)$ be a one parameter subgroup.
Then $\exp t\.x$ is grouplike for all $t$, i.e.,

\begin{align*}\exp t\.c(x)&=\exp c(t\.x)=c(\exp t\.x)=
 (\exp t\.x)\til\otimes(\exp t\.x)\\
&=\big((\exp t\.x)\til\otimes1\big)\big(1\til\otimes(\exp t\.x)\big)\\
&=\exp(t\.x\til\otimes1)\exp(1\til\otimes t\.x)
=\exp((t\.x\til\otimes1)+(1\til\otimes t\.x)) \\
&=\exp t\.(x\til\otimes 1+1\til\otimes x),\hbox{\quad for all $t\in\R$},
\end{align*}
\noindent
and this implies 
$$c(x)=x\til\otimes 1 + 1\til\otimes x$$ which  means $x\in P(A)$.

\ssk
If the weakly complete Hopf algebra $A$ arises as the dual of an
(abstract) Hopf algebra $H$ (i.e., a group object in $\mathcal{V}$), then
the members of $G(A)$ are multiplicative linear functionals on 
$H$, the so-called {\it characters} of $H$. Thus, Theorem \ref{7.8}
may be interpreted as saying that {\it the character group
of a Hopf algebra is a pro-Lie group.} (See also Theorem 5.6 in \cite{bog}.)

\ssk 
The simplest example is $A=\R[[x]]$, the algebra of formal power series
in one variable. As a vector space, $A$ is isomorphic to $\R^{\{1,x,x^2,\dots\}}$,
and this is weakly complete. Then $A\til\otimes A$ is isomorphic
to $\R[[y,z]]$, the formal power series algebra in two commuting variables $y$
and $z$.  The algebra morphism $\R[[x]]\to\R[[y,z]]$ generated by
   $x\mapsto y+z$ gives a
comultiplication $c\colon A\to A\til\otimes A$ making $A$ into 
a weakly complete Hopf algebra. The multiplicative subgroup $G(A)=
\{\exp t\.x: t\in\R\}$ is a Lie group isomorphic to $(\R,+)$ and
$P(A)\cong\R\.x$ is (trivially) a Lie algebra mapped by $\exp$
onto $G(A)$.

\msk
\noindent\textbf{Question 6.} (i) Is there a more elaborate duality theory of real Hopf algebras
and weakly complete Hopf algebras in which these facts on pro-Lie group and pro-Lie
algebra theory play a role?  \hfill\break
    (ii) Does the existence of weakly complete enveloping algebras
    of weakly complete Lie algebras secure for each pro-Lie algebra
    $L$ an associative weakly complete Hopf Algebra $U(L)$  such that
    $L$ is isomorphic to a closed Lie subalgebra of $P(U(L))$?

\section{ Further open questions}

\noindent In this last section we record a few additional  questions
on pro-Lie groups which do not fit naturally with
the material previously discussed in this paper, but are of some significance
to pro-Lie group theory. In so doing we rely on definitions and concepts
defined and discussed in \cite{PROBOOK}. 

\bsk
\noindent\textbf{Question 7.} Is an abelian prodiscrete compactly generated
group without nondegenerate compact subgroups a discrete group? 

\ssk
For a a definition and discussion of compactly generated groups,
see \cite{PROBOOK}, Definition 5.6, pp. 218ff. 
For a discussion of abelian compactly generated pro-Lie groups,
see \cite{PROBOOK}, Theorem 5.32, p. 236.

\msk
\noindent\textbf{Question 8.}  Is a compactly generated
 abelian prodiscrete compact--free 
group a finitely generated free abelian group?

\ssk
Note that by compact-free we mean the group has no nontrivial 
compact subgroups.

See \cite{PROBOOK}, Theorem 5.32,
the Compact Generation Theorem for Abelian Pro-Lie Groups, p.~236.

\ssk
In the proof of the structure of reductive pro-Lie algebras in  Theorem \ref{5.5}
(\cite{PROBOOK}, Theorem 7.27) one uses the lemma that in every finite-dimensional real semisimple Lie algebra every element is a sum of at
most two Lie brackets. 
 For brackets in semisimple Lie algebras see
\cite{PROBOOK}, Appendix 3, p.~651ff.

\msk\noindent\textbf{Question 9.}  Is every element in an arbitrary  
real semisimple Lie algebra 
a bracket?

\ssk
For the concept of  transfinitely solvable pro-Lie algebras and
pro-Lie groups, see \cite{PROBOOK}, Definition 7.32, pp. 285ff., 
respectively, pp.~420ff.
For the concept of  transfinitely nilpotent pro-Lie algebras and
pro-Lie groups, see \cite{PROBOOK}, pp. 296ff., respectively, pp.~443
ff.

\msk
\noindent\textbf{Question 10.}
 Is a transfinitely nilpotent connected pro-Lie group
pronilpotent?

\ssk
In Question 10 such a group  has to be prosolvable
since it is transfinitely solvable and then the
Equivalence Theorem for Solvability of Connected Pro-Lie Groups
10.18 of \cite{PROBOOK} applies. The impediment to a proof
is the failure of transfinite nilpotency to
be preserved by passing to quotients. Free topological groups
are free groups in the algebraic sense and thus 
are countably nilpotent;  but {\it every} topological 
group is a quotient of a free topological group and thus
of a transfinitely countably nilpotent topological group.

\ssk
For the definition of an analytic subgroup see \cite{PROBOOK},
 Definition 9.5 on p.~360.

\msk
 \noindent\textbf{Question 11.} Let $\h$ be a closed subalgebra 
of the Lie algebra $\g$ of a
connected pro-Lie group $G$. Let $A(\h)$ denote the analytic group
generated by $\h$. Is $\overline{A(\h)}/A(\h)$  abelian?

\ssk
(See \cite{PROBOOK}, Theorem 9.32.)

\msk\noindent\textbf{Question 12.}  Is there a satisfactory  theory of Polish 
 pro-Lie groups (respectively, separable  pro-Lie groups, or compactly generated
 pro-Lie groups), notably, in the connected case?

\msk

For information in \cite{PROBOOK} on the abelian case, 
see pp.~235ff.


\acknowledgments{Acknowledgements} We gratefully acknowledge helpful and
constructive
conversations we had with 
and  {\sc Saak S. Gabriyelyan} (Be'er Sheva, Israel) and {\sc Rafael Dahmen} (Darmstadt, Germany)  who also provided valuable assistance
in locating typographical errors in earlier versions of our
manuscript. We also thank the referees for finding typographical errors. Finally, the second author thanks the Technische Universit\"at Darmstadt for hospitality during which time the first draft of this survey article was written.





\break
\conflictofinterests{Conflicts of Interest}

The authors declare no conflict of interest.

\bibliographystyle{mdpi}
\makeatletter
\renewcommand\@biblabel[1]{#1. }
\makeatother

\end{document}